# Limits of Recursive Triangle and Polygon Tunnels

Florentin Smarandache
University of New Mexico, Gallup Campus, USA

**Abstract.**
In this paper we present unsolved problems that involve infinite tunnels of recursive triangles or recursive polygons, either in a decreasing or in an increasing way. The "nedians or order i in a triangle" are generalized to "nedians of ratio r" and "nedians of angle α" or "nedians at angle β", and afterwards one considers their corresponding "nedian triangles" and "nedian polygons".
This tunneling idea came from physics. Further research would be to construct similar tunnel of 3-D solids (and generally tunnels of n-D solids).

**A) Open Question 1 (Decreasing Tunnel).**

1. Let $\triangle ABC$ be a triangle and let $\triangle A_1B_1C_1$ be its **orthic triangle** (i.e. the triangle formed by the feet of its altitudes) and $H_1$ its <u>orthocenter</u> (the point on intersection of its altitudes).
   Then, let's consider the triangle $\triangle A_2B_2C_2$, which is the orthic triangle of triangle $\triangle A_1B_1C_1$, and $H_2$ its orthocenter.
   And the recursive tunneling process continues in the same way.
   Therefore, let's consider the triangle $\triangle A_nB_nC_n$, which is the orthic triangle of triangle $\triangle A_{n-1}B_{n-1}C_{n-1}$, and $H_n$ its orthocenter.
   a) What is the locus of the orthocenter points $H_1$, $H_2$, …, $H_n$, … ? {Locus means the set of all points satisfying some condition.}
   b) Is this limit:
   $$\lim_{n \to \infty} \triangle A_nB_nC_n$$
   convergent to a point? If so, what is this point?
   c) Calculate the sequences
   $$\alpha_n = \frac{area(\triangle A_nB_nC_n)}{area(\triangle A_{n-1}B_{n-1}C_{n-1})} \text{ and } \beta_n = \frac{perimeter(\triangle A_nB_nC_n)}{perimeter(\triangle A_{n-1}B_{n-1}C_{n-1})}.$$
   d) We generalize the problem from triangles to polygons. Let $AB…M$ be a polygon with $m \geq 4$ sides. From $A$ we draw a perpendicular on the next polygon's side $BC$, and note its intersection with this side by $A_1$. And so on. We get another polygon $A_1B_1…M_1$.
   We continue the recursive construction of this tunnel of polygons and we get the polygon sequence $A_nB_n…M_n$.
   d1) Calculate the limit:

$$\lim_{n \to \infty} \Delta A_n B_n ... M_n$$

d2) And the ratios of areas and perimeters as in question *c)*.
e) A version of this polygonal extension *d)* would be to draw a perpendicular from *A* not necessarily on the next polygon's side, but on another side (say, for example, on the third polygon's side) – and keep a similar procedure for the next perpendiculars from all polygon vertices *B, C*, etc.

In order to tackle the problem in a easier way, one can start by firstly studying particular initial triangles $\Delta ABC$, such as the equilateral and then the isosceles.

### B) Open Question 2 (Decreasing Tunnel).

2. Same problem as in Open Question 1, but replacing "orthic triangle" by "medial triangle", and respectively "orthocenter" by "center of mass (geometric centroid)", and "altitude" by "median". Therefore:

   Let $\Delta ABC$ be a triangle and let $\Delta A_1 B_1 C_1$ be its **medial triangle** (i.e. the triangle formed by the feet of its medians on the opposite sides of the triangle $\Delta ABC$) and $H_1$ its <u>center of mass</u> (or <u>geometric centroid</u>) (the point on intersection of its medians).
   Then, let's consider the triangle $\Delta A_2 B_2 C_2$, which is the medial triangle of triangle $\Delta A_1 B_1 C_1$, and $H_2$ its center of mass.
   And the recursive tunneling process continues in the same way.
   Therefore, let's consider the triangle $\Delta A_n B_n C_n$, which is the medial triangle of triangle $\Delta A_{n-1} B_{n-1} C_{n-1}$, and $H_n$ its center of massy.
   a) What is the locus of the center of mass points $H_1, H_2, ..., H_n, ...$ ?
      {This has an easy answer; all $H_i$ will coincide with $H_1$ (FS, IP).}
   b) Is this limit:
      $$\lim_{n \to \infty} \Delta A_n B_n C_n$$
      convergent to a point? If so, what is this point?
      {Same response; the limit is equal to $H_1$ (FS, IP).}
   c) Calculate the sequences
      $$\alpha_n = \frac{area(\Delta A_n B_n C_n)}{area(\Delta A_{n-1} B_{n-1} C_{n-1})} \text{ and } \beta_n = \frac{perimeter(\Delta A_n B_n C_n)}{perimeter(\Delta A_{n-1} B_{n-1} C_{n-1})}.$$
   d) We generalize the problem from triangles to polygons. Let *AB...M* be a polygon with $m \geq 4$ sides. From *A* we draw a line that connects *A* with the midpoint of *BC*, and note its intersection with this side by $A_1$. And so on. We get another polygon $A_1 B_1 ... M_1$.
      We continue the recursive construction of this tunnel of polygons and we get the polygon sequence $A_n B_n ... M_n$.

d1) Calculate the limit:
$$\lim_{n \to \infty} \Delta A_n B_n \ldots M_n.$$

d2) And the ratios of areas and perimeters of two consecutive polygons as in question *c)*.

e) A version of this polygonal extension *d)* would be to draw a line that connects *A* not necessarily on the midpoint of the next polygon's side, but with the midpoint of another side (say, for example, of the third polygon's side) – and keep a similar procedure for the next lines from all polygon vertices *B, C*, etc.

**C) Open Questions 3-7 (Decreasing Tunnels).**

3. Same problem as in Open Question 1, but considering a tunnel of **incentral triangles** and their incentral points, and their interior angles' bisectors.
   Incentral triangle is the triangle whose vertices are the intersections of the interior angle bisectors of the reference triangle $\Delta ABC$ with the respective opposite sides of $\Delta ABC$.
4. Same problem as in Open Question 1, but considering a tunnel of **contact triangles** (intouch triangles) and their incircle center points, and their interior angles' bisectors.
   A contact triangle is a triangle formed by the tangent points of the triangle sides to its incircle.
5. Same problem as in Open Question 1, but considering a tunnel of **pedal triangles** and a fixed point *P* in the plane of triangle $\Delta ABC$.
   A pedal triangle of P is formed by the feet of the perpendiculars from P to the sides of the triangle $\Delta ABC$.
6. Same problem as in Open Question 1, but considering a tunnel of **symmedial triangles.**
   "The **symmedial triangle** $\Delta K_A K_B K_C$ (a term coined by E.W. Weisstein [4]), is the triangle whose vertices are the intersection points of the symmedians with the reference triangle $\Delta ABC$."
7. Same problem as in Open Question 1, but considering a tunnel of **cyclocevian triangles.**
   A cyclocevian triangle of triangle $\Delta ABC$ with respect to the planar point *P* is the Cevian triangle of the cyclocevian conjugate of $P$.

**D) Open Questions 8-12 (Increasing Tunnels).**

8. Similar problem as in Open Question 1, but considering a tunnel of **anticevian triangles** of the triangle $\Delta ABC$ with respect to the same planar point *P*. For question *c)* and *d1)* only.
   The anticevian triangle of the given triangle $\Delta ABC$ with respect to the given point *P* is the triangle of which $\Delta ABC$ is the Cevian triangle with respect to *P*.
9. Similarly, but considering a tunnel of **tangential triangles.**
   The tangential triangle to the given triangle $\Delta ABC$ is a triangle formed by the tangents to the circumcircle of $\Delta ABC$ at its vertices.
10. Similarly, but considering a tunnel of **antipedal triangles.**

The antipedal triangle of the given triangle $\triangle ABC$ with respect to the given point $P$ is the triangle of which $\triangle ABC$ is the pedal triangle with respect to $P$.

11. Similarly, but considering a tunnel of **excentral triangles.**
    The excentral triangle (or tritangent triangle) of the triangle $\triangle ABC$ is the triangle with vertices corresponding to the excenters of $\triangle ABC$.

12. Similarly, but considering a tunnel of **anticomplementary triangles.**
    The anticomplementary (or antimedian) triangle of the triangle $\triangle ABC$ is the triangle formed by the parallels drew through the vertices of the triangle $\triangle ABC$ to the opposite sides.

**E) Open Questions Involving Nedians 13-14.**

a) One calls <u>nedians of order $i$</u> [see 4] of the triangle $\triangle ABC$ the lines that pass through each of the vertices of the triangle $\triangle ABC$ and divide the opposite site of the triangle into the ratio $i/n$, for $1 \leq i \leq n-1$.
   Let's generalize this to **nedians of ratio $r$**, which means lines that pass through each of the vertices of the given triangle $\triangle ABC$ and divide the opposite site of the triangle into the ratio $r$.
   We introduce the notion of *nedian triangles*, first the <u>interior nedian triangle of order $i$</u> (or more general <u>interior nedian triangle of ratio $r$</u>), which is the triangle formed by the three points of intersections of the three nedians of order $i$ (or respectively of the three nedians of ratio $r$), taken two by two;
   and that of <u>exterior nedian triangle of order $i$</u> (or more general <u>exterior nedian triangle of ratio $r$</u>), which is the triangle $\triangle A'B'C'$ such that $A' \in BC$, $B' \in CA$, and $C' \in AB$ - where $AA'$, $BB'$, and $CC'$ are nedians of order $i$ (respectively of ratio $r$) in the triangle $\triangle ABC$.

b) Another notion to introduce: **nedians of angle $\alpha$** (or **$\alpha$-nedians**), which are nedians that each of them forms the same angle $\alpha$ with its respective side of the triangle, i.e.
   $$\measuredangle(AA', AB) = \measuredangle(BB', BC) = \measuredangle(CC', CA) = \alpha.$$
   And associated with this we have <u>interior $\alpha$-nedian triangle</u> and <u>exterior $\alpha$-nedian triangle.</u>

c) And one more derivative to introduce now: <u>nedians at angle $\beta$ to the opposite side</u> (or <u>nedians-$\beta$</u>), which are of course nedians that form with the opposite side of the triangle $\triangle ABC$ the same angle $\beta$.
   {As a particular case we have the altitudes, which are nedians at an angle of $90°$ or $90°$-nedians.}
   And associated with this we have <u>interior nedian-$\beta$ triangle</u> and <u>exterior nedian-$\beta$ triangle.</u>

d) All these notions about nedians in a triangle can be extended to **nedians in a polygon**, and to the formation of corresponding **nedian polygons**.

Then:

13. Let $\triangle ABC$ be a triangle and let $\triangle A_1B_1C_1$ be its **interior nedian triangle of ratio r.**
    Then, let's consider the triangle $\triangle A_2B_2C_2$, which is the interior nedian triangle of order i of triangle $\triangle A_1B_1C_1$.
    And the recursive tunneling process continues in the same way.

Therefore, let's consider the triangle $\Delta A_nB_nC_n$, which is the interior nedian triangle of ordedr i of triangle the triangle $\Delta A_{n-1}B_{n-1}C_{n-1}$.
Same questions b)-e) as in Open Question 1.
14. Similar questions for **exterior nedian triangle of ratio r**.
15-16. Similar questions for **interior α-nedian triangle** and **exterior α-nedian triangle**.
16-17. Similar questions for **interior nedian-β triangle** and **exterior nedian-β triangle**.
18-23. Similar questions as the above 13-17 for the corresponding **nedian polygons**.

### F) More Open Questions.

The reader can exercise his or her research on other types of decreasing or increasing tunnels of special triangles (if their construction may work), such as the: extangential triangle, cotangential triangle, antisuplementary triangle, automedial triangle, altimedial triangle, circumpedal triangle, antiparalel triangle, Napoléon triangles, Vecten triangles, Sharygin triangles, Brocard triangles, Smarandache-Pătrașcu triangles (or orthohomological triangles[i]), Carnot triangle, Fuhrmann triangle, Morley triangle, Țiţeica triangle, Lucas triangle, Lionnet triangle, Schroeter triangle, Grebe triangle, etc.
{We don't present their definitions since the reader can find them in books of *Geometry of Triangle* or in mathematical encyclopedias, see for examples [1] and [6].}

### G) Construction.

Further research would be to construct similar tunnels of 3-D solids (and, more general, **tunnels of n-D solids** in $R^n$).

---

[i] We call two triangles, which are simultaneously orthological and homological, *orthohomological triangles (*or *Smarandache-Pătrașcu triangles* [2]*)*; for example: if the triangle $\Delta ABC$ is given and P is a point inside it such that its pedal triangle $\Delta A_1B_1C_1$ is homological with $\Delta ABC$, then we say that the triangles $\Delta ABC$ and $\Delta A_1B_1C_1$ are orthohomological.